\documentclass{article}
\usepackage{amsmath,amsfonts,amsthm,epsfig,amssymb,color}
\newtheorem{theorem}{Theorem}
\newtheorem{lemma}[theorem]{Lemma}

\newtheorem*{claim*}{Claim}
\newtheorem{corollary}[theorem]{Corollary}

\renewcommand\Pr{{\mathop{\mathbb P{}}\nolimits}}
\renewcommand\phi{\varphi}
\newcommand\eps{\varepsilon}
\newcommand\noproof{\hfill$\Box$}

\newcommand\col{\mathrm{col}}
\newcommand\pmax{p_{\max}}
\newcommand\VV{{\mathcal V}}
\newcommand\ob[1]{\partial^\infty{#1}}
\newcommand\ol[1]{\overline{#1}}

\newcommand\ceil[1]{\lceil #1 \rceil}
\newcommand\OS{O^*}
\newcommand\El{\mathrm{loc}}
\newcommand\Egl{E_{\mathrm{gl}}}
\newcommand\Eglp{E_{\mathrm{gl}}'}
\newcommand\la{\lambda}
\newcommand\cc{{\mathrm{c}}}
\newcommand\E{\operatorname{\mathbb E{}}}
\newcommand\vol{{\operatorname{\mathrm{vol}}}}
\newcommand\area{{\operatorname{\mathrm{area}}}}
\newcommand\PrTp{\Pr^{\TT(s)}_p}
\newcommand\PrTpp{\Pr^{\TT(s)}_{p'}}
\newcommand\PrTh{\Pr^{\TT(s)}_{1/2}}
\newcommand\Prp{\Pr_p}
\newcommand\Prh{\Pr_{1/2}}
\newcommand\Hb{H_{\mathrm b}}
\newcommand\Vb{V_{\mathrm b}}
\newcommand\Hw{H_{\mathrm w}}
\newcommand\Vw{V_{\mathrm w}}
\newcommand\Ep{\E_p}
\newcommand\Z{{\mathbb Z}}
\newcommand\RR{{\mathbb R}}
\newcommand\PP{{\mathcal P}}
\newcommand\PR{{\mathcal R}}
\newcommand\TT{{\mathbb T}}
\newcommand\pH{p_{\mathrm H}}
\newcommand\pT{p_{\mathrm T}}

\newcommand\bb[1]{\bigl(#1\bigr)}

\newcommand\Ed{d_2}
\newcommand\JMd{d_{\mathrm{JM}}}
\newcommand\JMn[1]{||#1||_{\mathrm{JM}}}
\newcommand\norm[1]{||#1||}
\newcommand\En[1]{||#1||_2}
\newcommand\whp{whp}
\newcommand\pbad{p_{\mathrm{bad}}}
\newcommand\pneut{p_{\mathrm{neutral}}}
\newcommand\pgood{p_{\mathrm{good}}}
\newcommand\db[1]{\bigl(#1\bigr)}
\newcommand\Ec{E_3^\mathrm{crude}}

\newcommand\DD{{\mathcal D}} 
\newcommand\pDD{{\mathcal D}'} 
\newcommand\PPop{\PP^+_{1,s-\delta}} 

\begin{document}
\title{Percolation on random Johnson--Mehl tessellations and related models}

\date{October 24, 2006; revised February 16 2007}

\author{B\'ela Bollob\'as%
\thanks{Department of Pure Mathematics and Mathematical Statistics,
Wilberforce Road, Cambridge CB3 0WB, UK}
\thanks{Department of Mathematical Sciences,
University of Memphis, Memphis TN 38152, USA}
\thanks{Research supported in part by NSF grants CCR-0225610 and DMS-0505550
and ARO grant W911NF-06-1-0076}
\and Oliver Riordan${}^*$%
\thanks{Research supported by a Royal Society Research Fellowship}
}
\maketitle

\begin{abstract}
We make use of the recent proof that the critical probability for
percolation on random Voronoi tessellations is $1/2$
to prove the corresponding result for random Johnson--Mehl tessellations, as
well as for two-dimensional slices of higher-dimensional Voronoi
tessellations. Surprisingly, the proof is a little simpler for these
more complicated models.
\end{abstract}

\section{Introduction and results}

The {\em Johnson--Mehl tessellation} of $\RR^d$ may be described as
follows: particles (nucleation centres) arrive at certain times
according to a spatial (deterministic or random) birth process on
$\RR^d$. When a particle arrives, it starts to grow a `crystal' at
a constant rate in all directions. Crystals grow only through
`vacant' space not yet occupied by other crystals; they stop growing
when they run into each other. Also, a new particle that arrives
inside an existing crystal never forms a crystal at all. This
generates a covering of $\RR^d$ by crystals meeting only in their
boundaries: every point of $\RR^d$ belongs to the crystal that first
reached it, or to the boundaries of two or more such crystals if it
is reached simultaneously by several crystals.

These tessellations were introduced by Johnson and Mehl \cite{JM39}
in 1939 as spatial models for the growth of crystals in metallic
systems. The same growth model (but not the resulting tessellation)
had been considered earlier by Kolmogorov~\cite{Kol37}; 
similar models were introduced independently by
Avrami \cite{Avr39}, \cite{Avr40}. Not
surprisingly, these models go under a variety of names (see the
references): in mathematics, they tend to be called
{\em Johnson--Mehl tessellations}, so this is the term we shall use here.
These models have been used to analyze a great variety of problems from
phase transition kinetics to polymers, ecological systems and DNA
replication (see Evans \cite{Eva93}, Fanfoni and Tomellini
\cite{FT98}, \cite{FT05}, Ramos, Rikvold and Novotny \cite{RRN99},
Tomellini, Fanfoni and Volpe \cite{TFV00}, \cite{TFV02}, and
Pacchiarotti, Fanfoni and Tomellini \cite{PFT05}, to mention only a
handful of papers);
mathematical properties of these tessellations have been studied
by Gilbert~\cite{Gilbert62}, Miles~\cite{Miles72},
M{\o}ller~\cite{Moller92,Moller96}, Chiu and Quine~\cite{CQ97}, and
Penrose~\cite{Penrose02}, among others.

If all particles arrive at the same time then we get a {\em Voronoi
tessellation}; this was introduced into crystallography by Meijering
\cite{Mei53} in 1953, although it had been studied much earlier by
Delesse \cite{Del1848}, Dirichlet \cite{Dir1850} and Voronoi
\cite{Vor1909}, in whose honour it is named. Random Voronoi
tessellations have been studied in numerous papers: for a list of
references, see \cite{BRbook};
here, we shall make heavy
use of the results in \cite{Voronoi}. For a discussion of many aspects of Voronoi
and related tessellations, including random Voronoi tessellations and Johnson--Mehl
tessellations, see the book by Okabe, Boots, Sugihara and Chiu~\cite{OBSC}.

In this paper we are mainly interested in {\em random Johnson--Mehl
tessellations of the plane}. As in almost all probabilistic models
in the literature, we shall assume that the birth process is a
time-homogeneous Poisson process of constant intensity, say,
intensity 1. Thus, our particles arrive randomly on the plane at
random times $t\ge 0$, according to a homogeneous Poisson process
$\PP$ on $\RR^2\times [0,\infty)$.

\begin{figure}[htb]
 \[\epsfig{file=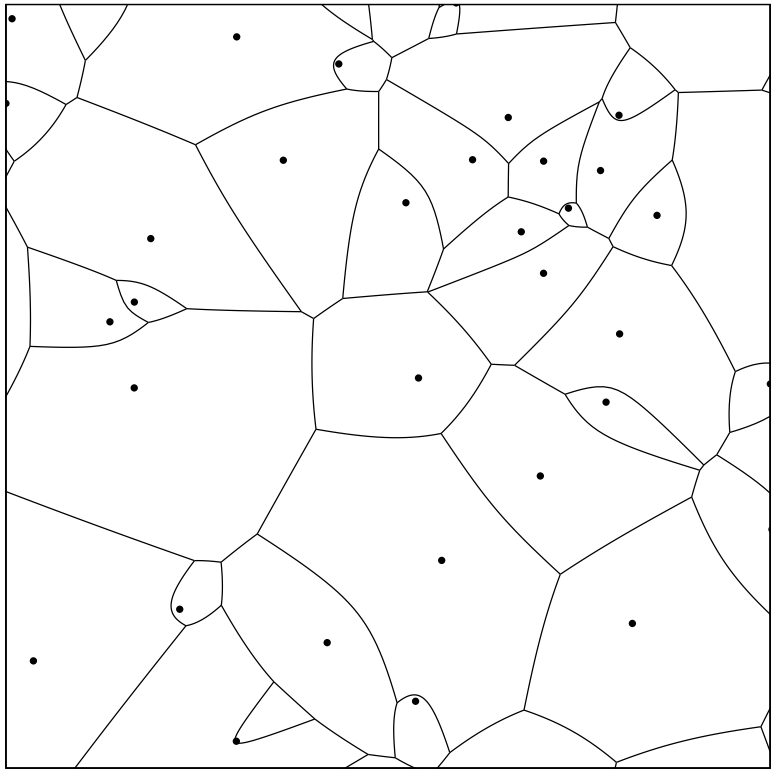,width=4in}\]
 \caption{Part of a random Johnson--Mehl tessellation of $\RR^2$. The dots
are the projections onto $\RR^2$ of those points $z$ of a Poisson process
in $\RR^2\times[0,\infty)$ for which the corresponding cell $V_z$ is non-empty.}\label{fig_JM}
\end{figure}

For each particle arriving at position $w\in \RR^2$
and time $t\ge 0$ we have a point $z=(w,t)\in \PP$.
The crystal associated to $z=(w,t)$ reaches a point $x\in \RR^2$ at time
$\Ed(x,w)+t$, where $\Ed$ denotes Euclidean distance.
(The subscript in the notation refers to the power in the norm,
not to the dimension: we shall write $d_p$ for the metric
on $\RR^d$ associated to the $\ell_p$-norm.)
Let $\JMn{\cdot}$ denote the norm on $\RR^3$ defined by
\[
 \JMn{(x_1,x_2,t)} = \sqrt{x_1^2+x_2^2}+|t| = \En{(x_1,x_2)} + |t|,
\]
and let $d=\JMd$ denote the corresponding distance.
Then the crystal $V_z$ associated to $z\in \PP$ may be written as
\begin{equation}\label{Vi}
 V_z = \bigl\{\ x\in \RR^2: d\db{(x,0),z} = \inf_{z'\in\PP} d\db{(x,0),z'} \ \bigr\}.
\end{equation}
A portion of a Johnson--Mehl tessellation is shown in Figure~\ref{fig_JM}.

In this way we see that the Johnson--Mehl tessellation of $\RR^2$
corresponds to a two-dimensional slice of the Voronoi tessellation
of $\RR^3$ defined with respect to a slightly unusual metric, $\JMd$;
only the absolute values of the time coordinates appear in
\eqref{Vi}, so up to rescaling (changing the density by a factor of
two), it makes no difference whether we take the particles to form a
Poisson process on $\RR^3$ or on $\RR^2\times [0,\infty)$.

In addition to studying Johnson--Mehl tessellations, we shall also
study two-dimensional slices of the usual random Voronoi
tessellation of $\RR^3$: the cells are defined exactly as in
\eqref{Vi}, but using the usual Euclidean metric $\Ed$ on $\RR^3$.
In fact, much of what we shall say will apply
to more general norms on $\RR^3$; however, the tessellation of
$\RR^2$ associated to a general norm on
$\RR^3$ is a rather unnatural object. Indeed, for a general norm,
the cells $V_z\subset \RR^2$ need not even be connected, and the
associated graph $G_{\PP}$ defined below need not be planar. For
this reason, we shall focus our attention on tessellations
associated to $\JMn{\cdot}$ and to the Euclidean norm. Another
example we shall consider is the norm $\ell_1$ on $\RR^3$, which may
be viewed as $\ell_1\oplus\ell_1$ on $\RR^2\oplus \RR^1$; the
associated tessellation of $\RR^2$ is a Johnson--Mehl type
tessellation in which crystals grow as squares whose side-lengths
increase at a constant rate.

As our main focus will be the Johnson--Mehl tessellation, we shall
always take $\PP$ to be a Poisson process on $\RR^2\times
[0,\infty)$, rather than on $\RR^3$, noting that the only effect on
the resulting tessellation of $\RR^2$ is a rescaling.

Having defined the cells $V_z$ associated to the points $z$ of a
Poisson process $\PP$, there is a natural way to construct an
associated graph $G_\PP$: the vertex set may be taken either to be the set of
$z\in\PP$ for which $V_z\ne\emptyset$, or all of $\PP$ (in which
case vertices corresponding to empty cells will be isolated). Two
vertices are adjacent if the corresponding cells meet, i.e., share one
or more boundary points. Ignoring probability zero events, as we may,
two vertices are adjacent if and only if
their cells have a common boundary arc. (Two cells may share more
than one boundary arc; there is an example in Figure~\ref{fig_JM}.)

Our aim is to study site percolation on the random graph $G_\PP$, or, equivalently,
`face percolation' on the tessellation $\{V_z\}$ itself. Let $0<p<1$ be a parameter.
We assign a {\em state}, {\em open} or {\em closed},
to each vertex of $G_\PP$, so that, given $\PP$, the states
of the vertices are independent, and each is open with probability $p$.
We are interested in the question `for which $p$ does $G_\PP$ contain an infinite
connected subgraph all of whose vertices are open?'. Equivalently, we may colour the cells
$V_z$ of the tessellation independently, taking each cell to be black with probability $p$ and white
otherwise, and we ask for which $p$ there is an unbounded {\em black component}.
We shall switch freely between these two viewpoints, writing $\Pr_p$ for the (common)
associated probability measure.

Let us say that a point $x\in \RR^2$ is {\em black} if it lies in a black cell,
and {\em white} if it lies in a white cell. Note that a point may be both black
and white, if it lies in the boundary of two cells.
Let $C_0$ be the set of points of $\RR^2$ joined to the origin by a
{\em black path}, i.e., a topological path in $\RR^2$ every point of which is black.
Let $z_0$ be the a.s. unique point of $\PP$ in whose cell the origin lies, and
let $C_0^G$ be the {\em open cluster} of $G_{\PP}$ containing $z_0$, i.e., the set
of all vertices of $G_{\PP}$ joined to $z_0$ by a path in the graph $G_{\PP}$
in which every vertex is open.
Since the cells $V_z$ are connected (see Section~\ref{sec_basics}),
the set $C_0\subset \RR^2$ is precisely
the union of the $V_z$ for $z\in C_0^G$.

Let
\[
 \theta(p)=\Prp(|C_0^G|=\infty) = \Prp(\,C_0\hbox{ is unbounded}\,),
\]
and let
\[
 \chi(p)=\Ep(|C_0^G|),
\]
where $\E_p$ is the expectation corresponding to $\Pr_p$.
Note that the graph $G_\PP$ depends on the metric $d$ as well as on $\PP$.
Thus $\theta(p)$ and $\chi(p)$ depend on $d$; most of the time, we suppress this dependence.

We say that our coloured random tessellation {\em percolates} if $\theta(p)>0$. It is easy
to see (from Kolmogorov's $0$-$1$ law, say) that,
in this case, the tessellation a.s. contains an unbounded
black component,
while if $\theta(p)=0$, then a.s. there is none.
We write
\[
 \pH=\pH(d)=\inf\{p:\theta(p)>0\}
\]
for the Hammersley critical probability associated to
percolation on our random tessellation, and
\[ 
 \pT=\pT(d)=\inf\{p:\chi(p)=\infty\}
\]
for the corresponding Temperley critical probability.

Our main aim in this paper is to determine the critical
probabilities for the Johnson--Mehl tessellation and for a
two-dimensional slice of the three-dimensional Voronoi tessellation. The
corresponding task for the random Voronoi tessellation associated to a
homogeneous Poisson process on $\RR^2$ was accomplished recently in
\cite{Voronoi}, where it was proved that $\pH=\pT=1/2$.

\begin{theorem}\label{th1}
Let $d$ denote either $\JMd$ or $\Ed$, let $\PP$ be a homogeneous Poisson process on $\RR^2\times[0,\infty)$
or on $\RR^3$, and
let $G_{\PP}$ be the graph associated to the tessellation $\{V_z\}$ of $\RR^2$ defined by \eqref{Vi}.
Then $p_H(d)=p_T(d)=1/2$.
More precisely, $\theta(p)>0$ if and only if $p>1/2$ and, for every $p<1/2$, there is a constant
$a=a(p)>0$ such that
\[
 \Prp\bb{\mathrm{size}(C_0)\ge n} \le \exp(-a(p)n)
\]
for all $n\ge 1$, where $\mathrm{size}(C_0)$ is the area of $C_0$, the diameter of $C_0$, or the
number $|C_0^G|$ of cells in $C_0$. 
\end{theorem}

The proof of the corresponding result for Voronoi tessellations
in~\cite{Voronoi} is rather lengthy. Much of this proof adapts
easily to the Johnson--Mehl setting, including, for example,
the analogue of the Russo--Seymour--Welsh Lemma. However, the hardest
part of the proof, a certain 
technical lemma, Theorem 6.1 in~\cite{Voronoi}, does not.
This result asserts that one can approximate
the continuous Poisson process $\PP$ by a suitable discrete process;
the proof of this extremely unsurprising statement makes up
a significant fraction of the length of~\cite{Voronoi}.
The analogue of this result for the Johnson--Mehl model is
Theorem~\ref{th_couplerobust} below; because the arguments
depend on the details of the geometry, a fresh proof
is required here. Surprisingly, although the Johnson--Mehl
model is more complicated that the Voronoi model, the proof
turns out to be simpler, though still not short. The key
difference is that we can use the third dimension of the model
to our advantage.

In the next section we describe basic properties of the Johnson--Mehl
model. In Section~\ref{sec_red} we outline the proof
of Theorem~\ref{th1}, assuming Theorem~\ref{th_couplerobust};
this part of the paper consists of a straightforward
adaptation of arguments from~\cite{Voronoi}. The heart of the
present paper is Section~\ref{sec_heart}, where we prove the technical
approximation lemma for the Johnson--Mehl model. In the final
section we discuss some generalizations.

\section{Basic properties}\label{sec_basics}

The probability that some point of the plane is equidistant from
four points of $\PP$ is zero. Hence, with probability $1$, at most three
cells $V_z$ of the tessellation associated to $\PP$ meet at any point.
We shall always assume that $\PP$ has this property.
Similarly, given any measure zero set $N$ (for example, the boundary
of a fixed rectangle), we may assume that no point of $N$
lies in three cells. Also, as any ball in $\RR^3$ contains only finitely
many points of $\PP$ (a.s. or always, depending on the definition of
a Poisson process one chooses), we shall assume that every
disk in $\RR^2$ meets finitely many cells $V_z$.

If we take our metric $d$ to be the Euclidean metric $\Ed$
or the $\ell_1$-metric $d_1$
(and take $\PP$ to be a Poisson process on $\RR^3$),
then the cells $V_z$ are two-dimensional sections of (bounded) convex
sets (in fact polyhedra) in $\RR^3$, and hence convex. For $d=\JMd$ this
is not true, but the cell $V_z$ associated to a point $z=(w,t)\in \PP$
is still a star domain, with centre $w$:
\begin{figure}[htb]
\centering
\input{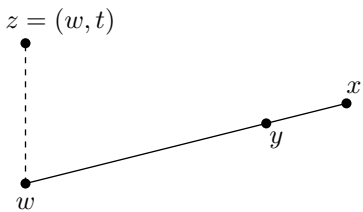}
\caption{A point $y$ on the line segment $wx$ in the plane. As we move
towards $w$ from $x$ at rate $1$, the $\JMd$-distance from $z$ decreases
at rate $1$. The $\JMd$-distance from any other $z'\in \PP$ decreases
at most this fast, so if $x$ lies in $V_z$ then so does~$y$.}\label{fig_star}
\end{figure}
if $x\in V_z$ and $y$ is a point on the line segment $wx$,
as in Figure~\ref{fig_star}, then
we have
\begin{multline*}
 \JMd\db{(y,0),z} = \En{y-w}+t = \En{x-w}-\En{y-x}+t \\
 = \JMd\db{(x,0),z} - \JMd\db{(x,0),(y,0)},
\end{multline*}
while for any $z'\in\RR^3$,
\[
 \JMd\db{(y,0),z'} \ge \JMd\db{(x,0),z'} - \JMd\db{(x,0),(y,0)}
\]
by the triangle inequality. Since $\JMd\db{(x,0),z}\le \JMd\db{(x,0),z'}$
for all $z'\in\PP$, the same inequality for $y$ follows, i.e., $y\in
V_z$. Thus $V_z$ is a star domain, and in particular $V_z$ is
connected.
Of course, the same argument applies to any metric $d$ on $\RR^3=\RR^2\oplus\RR$
that is the direct sum of a metric on $\RR^2$ and one on $\RR$.

\smallskip
Rather than first constructing a Poisson process $\PP$ on
$\RR^2\times[0,\infty)$, and then colouring the points
of $\PP$ black with probability $p$ and white with probability $1-p$,
equivalently we may start with two independent Poisson processes
$\PP^+$, $\PP^-$ with intensities $p$ and $(1-p)$, corresponding to
the black and white points, respectively. This is the viewpoint we
shall adopt most of the time.

In this viewpoint, our state space $\Omega$ consists of all pairs $(X^+,X^-)$
of discrete subsets of $\RR^2\times [0,\infty)$.
An event $E\subset \Omega$ is {\em black-increasing}, or simply
{\em increasing} if, whenever $(X_1^+,X_1^-)\in E$ and $(X_2^+,X_2^-)\in\Omega$
with $X_1^+\subset X_2^+$ and $X_1^-\supset X_2^-$, then $(X_2^+,X_2^-)\in E$.
In other words, $E$ is increasing if it is preserved by the addition of (black) points to $X^+$
and the deletion of (white) points from $X^-$.
If $x\in \RR^2$, then `$x$ is black' is an increasing event,
and so is any event of the form
`there exists a black path $P\subset \RR^2$ with certain properties.'

It is straightforward to check that Harris's Lemma concerning correlation of increasing
events extends to the present context; see~\cite{Voronoi}.

\begin{lemma}\label{l_Harris}
Let $E_1$ and $E_2$ be (black-)increasing events, and let $0<p<1$. Then
$\Prp(E_1\cap E_2) \ge \Prp(E_1)\Prp(E_2)$.
\noproof
\end{lemma}

Let us note the following simple fact for future reference.
\begin{lemma}\label{l_loc}
There is an absolute constant $A=A(d)$ with the following property:
let $S\subset \RR^2$ be a set with diameter at most $s$, and let
$\El(S)$ be the event that every point of $S$ is within $d$-distance $A(\log s)^{1/3}$
of some point of $\PP$. Then $\Pr(\El(S))=1-o(1)$ as $s\to\infty$.
\noproof
\end{lemma}

This lemma is a simple consequence of the basic properties of Poisson processes
(and is also a special case of a very weak form of a result of Penrose \cite{P1});
we omit the proof.

\section{Reduction to a coupling result}\label{sec_red}

In this section we present a proof of Theorem~\ref{th1}, assuming a
certain coupling result, Theorem~\ref{th_couplerobust} below, that
will allow us to discretize our Poisson process.  In a sense,
Theorem~\ref{th_couplerobust} is a technical lemma, and the arguments
in this section are the heart of the proof. However, as
in~\cite{Voronoi}, the hardest part of the overall proof is the proof
of Theorem~\ref{th_couplerobust}, presented in the next section.  The
arguments in this section are, {\em mutatis mutandis}, exactly the
same as those for random Voronoi percolation in~\cite{Voronoi},
so in places we shall only outline the details.

Given a rectangle $R=[a,b]\times [c,d]\subset \RR^2$, $a<b$, $c<d$,
let $H(R)=\Hb(R)$ be the event that there is a piecewise linear path $P\subset R$
joining the left- and right-hand sides of $R$ with every point of $P$ black.
When $\Hb(R)$ holds, we say that $R$ has a {\em black horizontal crossing}.
Let $V(R)=\Vb(R)$ be the event that $R$ has a {\em black vertical crossing}, defined
similarly. Also, let $\Hw(R)$ and $\Vw(R)$ denote the events that $R$ has a white horizontal
crossing or a white vertical crossing, respectively, defined in the obvious way.

Note that $H(R)=\Hb(R)$ is a black-increasing event. Also, from the topology of our tessellation,
$H(R)$ holds if and only if there is a sequence $z_1,\ldots,z_t$
of black points of $\PP$
such that the cells $V_{z_1}$ and $V_{z_t}$ meet the left- and right-hand sides of $R$,
respectively, and the cells of $V_{z_i}$ and $V_{z_{i+1}}$ meet at some point of $R$
for each $i$.

If no boundary point of $R$ lies in three or more Voronoi cells and no corner
of $R$ lies in two cells (which we may assume, as this event has probability $1$),
then from the topology of the plane exactly one of the events $\Hb(R)$ and $\Vw(R)$ holds,
so $\Prp(\Hb(R))+\Prp(\Vw(R))=1$.

Note that, from the symmetry of the model with respect to interchanging black and white,
$\Prp(\Vw(R))=\Pr_{1-p}(\Vb(R))$ for any $R$ and any $p$.
Furthermore, the metrics $d$ we consider are invariant under rotation (of the plane) through $\pi/2$,
so  $\Prp(\Hb(S))=\Prp(\Vb(S))$ for every square $S$. It follows that $\Prh(H(S))=1/2$.

Let $f_p(\rho,s)$ denote the $\Prp$-probability of the event $H([0,\rho s]\times [0,s])$,
i.e., the probability that a rectangle with aspect ratio $\rho$ and vertical side length
(or `scale') $s$ has a black horizontal crossing, and note that
\begin{equation}\label{s}
 f_{1/2}(1,s) = 1/2.
\end{equation}

The events $H(R)$ are defined in terms of the existence of certain black paths
in a certain black/white-colouring of the plane (in which some
points are both black and white). This random colouring has the following properties:
firstly, the event that any point (or given set of points)
is black is a black-increasing event,
so any two such events are positively correlated. Secondly, the distribution of the random
colouring is invariant under the symmetries of $\Z^2$, i.e., under translations (by integer
or in fact arbitrary vectors), under reflections in the axes, and under rotations through
multiples of $\pi/2$. Thirdly, well separated regions are asymptotically independent:
more precisely, let $\rho>0$ and $\eta>0$ be constants. Given $\eps>0$, if $s$ is large enough,
then for $R_1$ and $R_2$ two $\rho s$ by $s$ rectangles
separated by a distance of at least $\eta s$, and $E_1$ and $E_2$ any events determined by
the colours of the points (of $\RR^2$, not just of $\PP$)
within $R_1$ and $R_2$ respectively,
we have $|\Pr(E_1\cap E_2)-\Pr(E_1)\Pr(E_2)|\le \eps$.
To see this, note that when the event $\El(R_i)$ defined in Lemma~\ref{l_loc} holds,
the colouring of $R_i$ is determined by the positions and colours of the points of $\PP$
within distance $O((\log s)^{1/3})=o(s)$ of $R_i$.

As noted in \cite{Voronoi}, the properties above are all that is needed in the proof of Theorem 4.1
of that paper, which thus carries over to the present setting.
\begin{theorem}\label{th_RSW}
Let $0<p<1$ and $\rho>1$ be fixed. If $\liminf_{s\to\infty} f_p(1,s)>0$,
then $\limsup_{s\to\infty} f_p(\rho,s)>0$. \noproof
\end{theorem}
Together with \eqref{s}, Theorem~\ref{th_RSW} has the following corollary.
\begin{corollary}\label{c}
Let $\rho>1$ be fixed. There is a constant $c_0=c_0(\rho)>0$ such that
for every $s_0$ there is an $s>s_0$ with $f_{1/2}(\rho,s)\ge c_0$.
\noproof
\end{corollary}

As in the context of ordinary Voronoi percolation, to prove Theorem~\ref{th1}
it suffices to prove the following result, analogous to Theorem 7.1 of \cite{Voronoi}.
\begin{theorem}\label{thabove}
Let $\rho>1$, $p>1/2$, $c_1<1$ and $s_1$ be given.
There is an $s>s_1$ such that $f_p(\rho,s)>c_1$.
\end{theorem}
Theorem~\ref{th1} may be deduced from Theorem~\ref{thabove} by using the idea
of $1$-independent percolation. The argument is exactly the same as in the Voronoi setting,
so we shall not give it.

In the light of the comments above, our task is to deduce
Theorem~\ref{thabove} from Corollary~\ref{c}. The basic
idea is simple: for $s$ large, we shall show that a small increase in $p$
greatly increases $f_p(\rho,s)=\Prp(H(R))$, where $R$ is a $\rho s$ by $s$ rectangle.
If $H(R)$ were a symmetric event in a discrete product space, then this would
be immediate from the sharp-threshold result of Friedgut and Kalai~\cite{FK}.
Unfortunately, $H(R)$ is neither symmetric nor an event in a discrete product space,
so we have two difficulties to overcome. The first is easily dealt with, by working
on the torus.

Let $\TT(s)$ denote the $s$ by $s$ {\em torus}, i.e., the quotient of
$\RR^2$ by the equivalence relation $(x,y)\sim (x',y')$ if
$x-x',y-y'\in s\Z$. Instead of $\RR^2\times [0,\infty)$, we shall work
in the `thickened torus' $\TT(s)\times [0,t]$. Note that we do not
wrap around in the third direction. It turns out that
the precise thickness $t$ is not important in the arguments that follow:
we could use any thickness $t$ larger than a certain
constant times $(\log s)^{1/3}$ but bounded by a power of $s$. For simplicity
we shall set $t=s$, working in $\TT(s)\times [0,s]$ throughout.

Let us write $\PrTp$ for the probability measure associated to a Poisson process $\PP$ on $\TT(s)\times[0,s]$
of intensity $1$ in which each point is coloured black with probability $p$ and white otherwise,
independently of the process and of the other points. Alternatively, $\PrTp$ is the probability measure associated
to a pair $(\PP^+,\PP^-)$ of independent Poisson processes on $\TT(s)\times[0,s]$ with intensities
$p$ and $1-p$, respectively. Our metric $d=\JMd$ (or $\Ed$, or $d_1$)
induces a metric on $\TT(s)\times[0,s]\subset \TT(s)\times \RR$
in a natural way. Thus, associated to $\PrTp$ we have a random
black/white-coloured tessellation of $\TT(s)$ by
the Voronoi cells associated to $(\PP,d)$.

If we restrict our attention to a region that does not come close to `wrapping around' the torus,
then $\PrTp$ and $\Prp$ are essentially equivalent. More precisely, identifying $\TT(s)\times [0,s]$ with
$[0,s)^2\times [0,s]\subset \RR^3$, we may couple the measures $\PrTp$ and $\Prp$ by realizing
our coloured Poisson process on $\TT(s)\times[0,s]$
as a subset of that on $\RR^2\times[0,\infty)$. Let $\eps>0$ be fixed, and let $R=[\eps s,(1-\eps)s]^2$.
Whenever the event $\El(R)$ defined in Lemma~\ref{l_loc} holds, the colour of every point of $R$
is determined by the restriction of the Poisson process to $[0,s)^2\times [0,s]$, so the colourings
of $R$ associated to the measures $\PrTp$ and $\Prp$ coincide.
Hence, Lemma~\ref{l_loc} has the following consequence.

\begin{lemma}\label{RT}
Let $0<a,b<1$ be constant, and let $R_s$ be an $a s $ by $b s$ rectangle.
Then for every $p$ we have
\[
 \PrTp(H(R_s)) = \Prp(H(R_s)) +o(1)
\]
as $s\to\infty$.
\noproof
\end{lemma}

Lemma~\ref{RT} says that when studying crossings of rectangles,
we can work on the torus instead of in the plane. On the torus,
there is a natural way to convert $H(R)$ into a symmetric event; we shall return
to this shortly.

As in \cite{Voronoi}, we wish to apply a Friedgut--Kalai sharp-threshold result from~\cite{FK}.
A key step is to approximate our Poisson process $\PP$ on $\TT(s)\times [0,s]$ by a discrete
process. Given $\delta=\delta(s)>0$ with $s/\delta$ an integer, partition $\TT(s)\times [0,s]$
into $(s/\delta)^3$ cubes $Q_i$ of side-length $\delta$ in the natural way. (We may ignore the boundaries
of the cubes, since the probability that $\PP$ contains a point in any of these boundaries is $0$.)
As in \cite{Voronoi}, the {\em crude state} of a cube $Q_i$ is {\em bad} if $Q_i$ contains
one or more points of $\PP^-$, {\em neutral} if $Q_i$ contains no points of $\PP^-\cup\PP^+=\PP$,
and {\em good} if $Q_i$ contains one or more points of $\PP^+$ but no points of $\PP^-$.
Let $\delta=\delta(s)$ be a function of $s$ that tends
to $0$ as $s\to\infty$ (later, $\delta(s)$ will be a small negative power of $s$); all asymptotic
notation refers to the $s\to\infty$ limit. Writing $\gamma=\delta^3$,
since $\delta(s)\to 0$, 
each $Q_i$ is bad, neutral or good with respective probabilities
\begin{eqnarray}
 \pbad &=& 1-\exp\big({-}\gamma(1-p)\big) \sim \gamma(1-p), \nonumber \\
 \pneut &=& \exp(-\gamma),\label{pbng} \\
 \pgood &=& \exp\big({-}\gamma(1-p)\big)\big(1-\exp(-\gamma p)\big) \sim \gamma p.\nonumber
\end{eqnarray}
Also, the crude states of the $Q_i$ are independent.

Writing $N=(s/\delta)^3$ for the number of cubes, and representing bad, neutral and
good states by $-1$, $0$ and $1$ respectively,
the measure $\PrTp$ induces a product measure on the set $\Omega_N=\{-1,0,1\}^N$ of crude states.

To prove results about the continuous process, we shall pass to the discrete setting and then back;
starting from a realization $(\PP_1^+,\PP_1^-)$ of our Poisson process, first
we generate the corresponding
crude states, and then we return to a possibly different realization $(\PP_2^+,\PP_2^-)$ consistent with
the same crude states. An event
such as $H(R)$ need not survive these transitions: a point $x$ or path $P$ may be black with respect
to $(\PP_1^+,\PP_1^-)$ but not with respect to $(\PP_2^+,\PP_2^-)$. To deal with this problem, we
consider a `robust' version of the event that a point or path is black.

Given $\eta>0$, let us say
that a point $x\in \TT(s)$ is {\em $\eta$-robustly black}
with respect to $(\PP^+,\PP^-)$ if the closest
point of $\PP^+$ to $x$ is at least a distance $\eta$ closer than the closest point of $\PP^-$,
where all distances are measured in the metric $d$.
[Note that whenever $x$ is $\eta$-robustly black, the entire $(\eta/2)$-neighbourhood
of $x$ is black. There is no reverse implication: for any $\eta>0$ and any $r$,
it is possible for the $r$-neighbourhood of $x$ to be black without
$x$ being $\eta$-robustly black.]
A path $P$ is {\em $\eta$-robustly black} if every point of $P$ is $\eta$-robustly black.
Set
\[
 C_d = \sup\{ d(x,y) : x,y\in [0,1]^3\} < \infty.
\]
If $(\PP_i^+,\PP_i^-)$, $i=1,2$, are realizations of our Poisson process on $\TT(s)\times[0,s]$ consistent with the
same crude state, and a point $x\in \TT(s)$ is $(2C_d\delta)$-robustly black with respect to $(\PP_1^+,\PP_1^-)$,
then it is easy to check that $x$ is black with respect to $(\PP_2^+,\PP_2^-)$.

Our starting point for the proof of Theorem~\ref{thabove} is Corollary~\ref{c}, which gives
us (with reasonable probability) a certain black path. Fortunately, we can `bump up' a black
path to a robustly black path at the cost of increasing $p$ slightly, using the following
analogue of Theorem 6.1 of \cite{Voronoi}. Here, and in what follows, we say
that an event holds {\em with high probability}, or {\em whp}, if it has
probability $1-o(1)$ as $s\to\infty$ with any other parameters fixed.

\begin{theorem}\label{th_couplerobust}
Let $d$ denote either $\JMd$ or $\Ed$, and
let $0<p_1<p_2<1$ and $\eps>0$ be given. Let $\delta=\delta(s)$ be
any function with $0<\delta(s)\le s^{-\eps}$.  We may construct in
the same probability space Poisson processes $\PP_1^+$, $\PP_1^-$,
$\PP_2^+$ and $\PP_2^-$ on $\TT(s)\times [0,s]$ of intensities $p_1$, $1-p_1$, $p_2$ and $1-p_2$,
respectively, so that $\PP_i^+$ and $\PP_i^-$ are independent for $i=1,2$,
and the following global event $\Egl$ holds \whp\ as $s\to\infty$: for
every piecewise-linear path $P_1\subset\TT(s)$ which is black with respect to
$(\PP_1^+,\PP_1^-)$ there is a piecewise-linear path $P_2\subset\TT(s)$ which is
$(2C_d\delta)$-robustly black with respect to $(\PP_2^+,\PP_2^-)$, such
that every point of $P_2$ is within distance $\log s$ of some
point of $P_1$ and {\em vice versa}.
\end{theorem}

The proof of this result is a little involved, and will be given in the next section.
This is the only part of the present paper that is essentially different from
the arguments for usual Voronoi tessellations given in~\cite{Voronoi}.

Before turning to the proof of Theorem~\ref{th_couplerobust}, let us outline
how Theorem~\ref{thabove} follows.
The proof is exactly the same as that in Section 7 of~\cite{Voronoi},
{\em mutatis mutandis}: we use Corollary~\ref{c} and
Theorem~\ref{th_couplerobust} in place of their analogues Corollary
4.2 and Theorem 6.1 of~\cite{Voronoi}; we write $(2C_d\delta)$-robustly
black in place of $4\delta$-robustly black; $\gamma=\delta^3$,
the volume of each small cube $Q_i$, replaces $\gamma=\delta^2$,
the volume of a small square $S_i$ in~\cite{Voronoi}; finally,
$N=(s/\delta)^3$, the number of cubes $Q_i$, replaces $N=(s/\delta)^2$,
the number of squares $S_i$.

Very roughly, the strategy of the proof is as follows (for details
see~\cite{Voronoi}). Fix $p>1/2$, let $\eps>0$ be chosen below
(depending on $p$), and let $s$ be `sufficiently large'.
From Corollary~\ref{c}, after increasing $s$, if necessary, the crossing probability
$f_{1/2}(10,s/13)$ is at least some positive absolute constant.
By Lemma~\ref{RT}, it follows that in the torus, i.e., in the measure
$\PrTh$, the probability that a given $10s/13$ by $s/13$ rectangle
has a black horizontal crossing is also at least a positive constant.
Set $p'=(p+1/2)/2$, say, so that $1/2<p'<p$,
and set $\delta=s^{-\eps}$, decreasing $\delta$ slightly if necessary so that
$s/\delta$ is an integer.
Using Theorem~\ref{th_couplerobust}, we can convert a black
path in $\PrTh$ to a `nearby' robustly black path in $\PrTpp$:
it follows that the $\PrTpp$-probability that a given $3s/4$ by $s/12$
rectangle has a $(2C_d\delta)$-robustly black horizontal crossing
is not too small, i.e., is at least some constant $c>0$.

The sharp-threshold result we shall need is Theorem 2.2 of~\cite{Voronoi},
a simple modification of a result of Friedgut and Kalai, Theorem 3.2 of~\cite{FK}.
Consider the state space $\Omega_N=\{-1,0,1\}^N$ with a product measure,
in which the coordinates are independent and identically distributed.
An event $E$ in this space is {\em increasing} if $\omega=(\omega_i)_{i=1}^n\in E$
and $\omega_i\le \omega_i'$ for every $i$ imply $\omega'\in E$. Also, $E$ is {\em symmetric}
if there is a group acting transitively on the coordinates $1,2,\ldots,N$
whose induced action on $\Omega_N$ preserves $E$. 
Roughly speaking, Theorem 2.2 of~\cite{Voronoi}
says that if $E$ is a symmetric increasing event in $\Omega_N$,
and we consider product measures on $\Omega_N$ in
which the probability that a given coordinate is non-zero is `small',
say bounded by $\pmax$, then increasing the probability that each coordinate is $1$
by at least $\Delta$ and decreasing the probability that it is $-1$ by at least
$\Delta$ is enough to increase the probability of $E$ from $\eta$ to $1-\eta$,
where
\[
 \Delta= C\log(1/\eta) \pmax\log(1/\pmax)/\log N
\]
and $C$ is constant. (For details, see~\cite{Voronoi}.)

To apply the result above, we need a symmetric
event in a discrete product space. To achieve symmetry, following the notation
in Section 7 of~\cite{Voronoi}, we simply consider
the event $E_3$ that {\em some} $3s/4$ by $s/12$ rectangle in $\TT(s)$ has a
robustly black horizontal crossing. To convert to a discrete product
space, we divide $\TT(s)\times[0,s]$ into $N=(s/\delta)^3$ cubes of
volume $\gamma=\delta^3$, and consider the crude state of each cube as
defined above.  Let $\Ec$ be the event that the crude states of the cubes
are consistent with $E_3$, which may be naturally identified
with an event in $\Omega_N$.
Note that $\PrTpp(\Ec)\ge \PrTpp(E_3)\ge c>0$.

The Friedgut--Kalai result
implies that a small increase in the probability of black points
increases the probability of $\Ec$
dramatically (details below).  It follows that $\PrTp(\Ec)$ is very
close to $1$ if $s$ is large. Hence, there is a
very high $\PrTp$-probability that some $3s/4$ by $s/12$ rectangle has a black
crossing.  (Not necessarily a robustly black crossing: in passing to
the discrete approximation and back again, points of our Poisson process may
move slightly. However, as noted above, any crossing that was robustly black
remains black.) By a simple application of the square-root trick, one
can deduce that the $\PrTp$-probability that a {\em fixed} $s/2$ by
$s/6$ rectangle in $\TT(s)$ has a black horizontal crossing is also
very close to $1$.  Finally, using Lemma~\ref{RT} again it follows
that $f_p(3,s/6)$ can be made arbitrarily close to $1$,
and Theorem~\ref{thabove} follows.

Turning to the quantitative application of the sharp-threshold result,
we may take $\eta$ to be a (very small) absolute constant.
Each cube $Q_i$ is very small, and the probability
that a cube $Q_i$ is either good or bad is at most $\gamma=\delta^3$, so we may take $\pmax=\gamma$.
Also, passing from $\PrTpp$ to $\PrTp$ increases the probability that a given cube
is good, and decreases the probability that it is bad, by roughly $(p-p')\gamma$;
see \eqref{pbng}.
Hence, to deduce Theorem~\ref{thabove} from the Friedgut--Kalai result,
we need
\[
 (p-p')\gamma \ge C' \gamma \log(1/\gamma)/\log N,
\]
for some constant $C'$. With $p$ and $p'$ fixed, this reduces to
$C''\log(1/\gamma)/\log N<1$. Since $N=s^{\Theta(1)}$ and $\gamma=s^{-\Theta(\epsilon)}$,
this condition can be met by choosing $\eps$ sufficiently small.
Note that it is irrelevant whether $\gamma=\delta^3$ and $N=(s/\delta)^3$,
as here, or $\gamma=\delta^2$ and $N=(s/\delta)^2$, as in~\cite{Voronoi}. Indeed,
this part of the argument works unchanged in any dimension: the key point
is that we can afford only to discretize to a scale $\delta$ given by an
arbitrarily small negative power of $s$. Fortunately, Theorem~\ref{th_couplerobust}
applies for such a $\delta$.

\section{Replacing black paths by robustly black paths}\label{sec_heart}

It remains only to prove Theorem~\ref{th_couplerobust}. Roughly speaking,
this states that a small increase in the probability $p$ that each point is black
allows any black path to be replaced by a nearby robustly black path.
The proof, to which this section is devoted, turns
out to be the hardest part of the paper. 

We shall use the following fact about random Voronoi
tessellations in three dimensions; here $O$ denotes the origin.

\begin{theorem}\label{th_faces}
Let $\PP$ be a homogeneous Poisson process on $\RR^3$ of intensity
$1$, and let $d$ denote $\Ed$ or $\JMd$.  For $A>0$ let $E_k=E_{k,A}$ be the event that
$\PP$ contains $k$ points $P_1,\ldots,P_k$ with the following
property: there are points $Q_1,\ldots,Q_k\in \RR^3$ and real numbers
$0<r_1,\ldots,r_k<Ak^{1/3}$ such that $d(O,Q_i)=d(P_i,Q_i)=r_i$ for every
$i$, and $d(P_j,Q_i)\ge r_i$ for all $i$ and $j$.  If $A>0$ and $C>0$ are constant,
then
\[
 \Pr(E_k) =o(e^{-Ck})
\]
as $k\to\infty$.
\end{theorem}
Theorem~\ref{th_faces} is essentially equivalent to the following 
statement:
if $V_0$ is the cell of the origin in the Voronoi tessellation of $\RR^3$ defined 
using the point set
$\PP\cup\{O\}$ and metric $d$, then the probability that $V_0$ has at
least $k$ faces is $o(e^{-Ck})$, for any constant $C$.
In fact, Theorem~\ref{th_faces} easily implies this:
suppose that $V_0$ has $k$ faces, let $P_1,\ldots,P_k$ be the centres of the adjacent
Voronoi cells $V_1,\ldots,V_k$, and let $Q_i$ be a common point of $V_0$ and $V_i$ for each $i$.
Then $O$ and $P_i$ are the two closest points of $\PP\cup\{O\}$
to $Q_i$, so either $E_k$ holds,
or we have $d(O,Q_i)=d(P_i,Q_i)>Ak^{1/3}$ for some $i$. But in the latter
case there is a  ball $B$ (defined with respect to $d$) of radius $r=Ak^{1/3}$ meeting the origin
and containing no points of $\PP$. Placing $O(1)$ balls $B_i$ of radius $r/2$ so that any such ball $B$
contains one of the $B_i$, the probability that some $B_i$ contains no points of $\PP$
is $O\bb{\exp(-\vol(B_i))}$, which is much smaller than $e^{-Ck}$ if we choose $A$ large enough.

Let us remark that the statement above may well be known, at least for $d=\Ed$.
In two dimensions, very precise results are known; Hilhorst~\cite{Hilhorst}
has shown that the probability $p_k$
that the cell of the origin has exactly $k$ faces satisfies
\[
 p_k = \frac{A}{4\pi^2} \frac{(8\pi^2)^k}{(2k)!} \bb{1+O(k^{-1/2})}
\]
as $k\to\infty$, where $A$ is a certain constant given as an infinite product.
It is very likely that similar $(c/k)^k$ asymptotics hold
for Johnson--Mehl tessellations; to keep the proof simple, we prove
only the much weaker $o(e^{-Ck})$ bound.

\begin{proof}[Proof of Theorem~\ref{th_faces}]
The strategy of the proof is as follows: we shall show that if $E_k$ holds, and certain `bad'
events $B_1,\ldots,B_4$ of probability $o(e^{-Ck})$ do not hold,
then an impossible situation arises. In the proof, by
`distance' we mean the usual Euclidean distance. Constants in $O(\cdot)$ notation may depend on $A$ and
$C$ unless explicitly stated otherwise.

For $a>0$ constant, let $N_a$ be the number of points $z$ of $\PP$ with
$d(O,z)<ak^{1/3}$. Then $\E(N_a)=\Theta(a^3k)$, so, choosing $a$ small
enough, $\E(N_a)\le k/8$. Fixing such an $a$ from now on,
let $B_1$ be the `bad' event
\[
 B_1 = \{N_a\ge k/4\}.
\]
The
probability that a Poisson random variable with mean $k/8$ exceeds $k/4$
is $k^{-\Theta(k)}=o(e^{-Ck})$, so $\Pr(B_1)=o(e^{-Ck})$.

Let $S$ be a (Euclidean) ball of radius $\Theta(k^{1/3})$ containing $\{z:d(O,z)\le 2Ak^{1/3}\}$,
noting that if $E_k$ holds, then all $P_i$ lie in $S$.
Let $N$ be the number of points of $\PP$ in $S$, so $N$ has a Poisson distribution
with mean $\vol(S)=\Theta(k)$.
Setting $C_1=2\vol(S)/k$, let $B_2$ be the event
\[
 B_2 = \{N\ge C_1k\}.
\]
Then
$\Pr(B_2)=k^{-\Theta(k)}=o(e^{-Ck})$.

Given $N=|\PP\cap S|$, we may generate $\PP\cap S$ as a sequence $z_1,\ldots,z_N$, where the
$z_i$ are independent and each $z_i$ is chosen uniformly from $S$.
Let $c$ be a small constant to be chosen below, and 
let $B_3$ be the event that there are at least $k/4$ points $z_i\in \PP\cap S$
with $\Ed(z_i,z_j)\le c$ for some $j<i$:
\[
 B_3 =\bigl\{\  | \{z_i: \exists j<i \hbox{ with }\Ed(z_i,z_j)\le c \} | \ge k/4 \ \bigr\}.
\]
Given $N$ and $z_1,\ldots,z_{i-1}$, the probability that $\Ed(z_i,z_j)\le c$ for some $j<i$
is at most $p_i=(i-1)4\pi c^3/(3\vol(S))$.
Conditional on $N$, the number of points $z_i$ within distance $c$ of an earlier $z_j$
is thus dominated by a binomial distribution with parameters $N$ and $p_N$, whose
mean $\mu_N$ is at most $N^2 4\pi c^3/(3\vol(S))=\Theta(N^2c^3/k)$.
Whenever $B_2$ does not hold, we have $N\le C_1 k$, and so $\mu_N=O(c^3k)$.
Hence, $\Pr(B_2^\cc\cap B_3)$ is at most the probability that a certain Binomial distribution
with mean $O(c^3k)$ exceeds $k/4$. Choosing $c$ small enough, this probability is $o(e^{-Ck})$,
so $\Pr(B_3)\le \Pr(B_2)+\Pr(B_2^\cc\cap B_3)=o(e^{-Ck})$.

Let $c_1$ be a small constant to be chosen below, and let $N'$ be the
number of points $P\in \PP\cap S$ for which the angle between $OP$ and
the plane is within $c_1$ of $\pi/4$. Then $N'$ has a Poisson distribution
with mean $\Theta(c_1\vol(S))=\Theta(c_1k)$. Choosing $c_1$ small enough,
$\E(N')\le k/8$. Let $B_4$ be the event
\[
 B_4 = \{ N'\ge k/4 \},
\]
noting that $\Pr(B_4)=o(e^{-Ck})$.

Suppose that $E_k$ holds and that none of the events $B_i$, $1\le i\le 4$, holds.
To complete the proof of Theorem~\ref{th_faces}, it suffices to deduce a contradiction.
As $B_1$, $B_3$ and $B_4$ do not hold, there is a subset $X$ of $\{P_1,\ldots,P_k\}$
of size at least $k-3k/4=k/4$ with the following properties: for every $P_i\in X$ we have
$d(O,P_i)\ge ak^{1/3}$, the angle $\theta_i$ between $OP_i$ and the plane
is at most $\pi/4-c_1$ or at least $\pi/4+c_1$, and no two $P_i$, $P_j\in X$ are
within distance $c$.

Projecting the $\Theta(k)$ points $P_i\in X$ from the origin onto the unit sphere centred
at the origin, we find two points $P_i$, $P_j\in X$ whose projections are within
distance $\Theta(k^{-1/2})$. Without loss of generality we may assume that $i=1$, $j=2$
and that $P_2$ is at least as far from $O$ as $P_1$ is.
Let $P_1'$ be the point of the line segment $OP_2$ at Euclidean distance $\Ed(O,P_1)$
from $O$; see Figure~\ref{fig_PPQ}.
\begin{figure}[htb]
\centering
\input{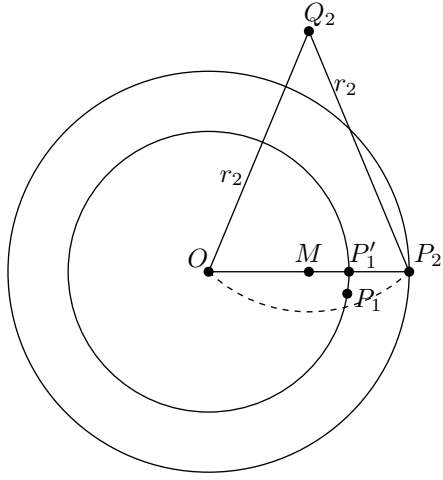}
\caption{Two points $P_1$, $P_2$ of $X$ such that the angle $P_1OP_2$ is very small.
$P_1'$ is the point of $OP_2$ at distance $d(O,P_1)$ from $O$. As $P_1OP_2$
is small, the distance $P_1P_1'$ is $o(1)$. Since $P_1P_2\ge c$,
it follows that $P_1'P_2\ge c/2$. The points $O$ and $P_2$ are equidistant
from $Q_2$; the dotted line is part of the sphere $d(Q_2,\cdot)=d(Q_2,O)$
(curvature exaggerated), and $M$ is the midpoint of $OP_2$.
Using convexity, we can show that $P_1$ lies inside the corresponding ball.}\label{fig_PPQ}
\end{figure}
Then $\Ed(P_1,P_1')=O(k^{-1/2}\Ed(O,P_1)) = O(k^{-1/2}k^{1/3})=o(1)$.
As $\Ed(P_i,P_j)\ge c$, for $k$ sufficiently large it follows that
$\Ed(P_1',P_2)\ge c/2$, i.e., that
$\Ed(O,P_2) \ge \Ed(O,P_1)+c/2$.

Recall that $O$ and $P_2$ are two points on the surface of a ball $B$ in the metric
$d$ with centre $Q_2$ and radius $r_2=\Theta(k^{1/3})$. We claim that the midpoint $M$
of $OP_2$ is `well inside' this ball, i.e., that
\begin{equation}\label{mi}
 d(M,Q_2)\le (1-\eps)r_2
\end{equation}
for some constant $\eps>0$ depending only on $d$, $A$ and $C$.
Recall that $d(O,P_2)\ge ak^{1/3}$ and $r_2\le Ak^{1/3}$,
so $d(O,P_2)=\Theta(r_2)$. If $d=\Ed$, then \eqref{mi} follows immediately: in fact,
it follows immediately for any metric $d$ defined by a `strictly convex' norm, i.e., one whose
unit sphere contains no line segments.

For $d=\JMd$, \eqref{mi} follows using the additional fact that the angle
between $OP_2$ and the plane is not within $c_1$ of $\pi/4$: the unit sphere
of the corresponding norm is a double cone. The only line segments this contains
are at an angle $\pi/4$ to the plane. Let $H$ be the set of
pairs $(A,B)$ of points on this unit sphere with $\Ed(A,B)$ at least some small
constant and the angle between $AB$ and the plane lying outside
the interval $(\pi/4-c_1,\pi/4+c_1)$.
Then $\JMn{(A+B)/2}<1$ for all $(A,B)\in H$, and hence, by compactness,
$\JMn{(A+B)/2}\le 1-\eps$ for all $(A,B)\in H$, for some $\eps>0$.
Relation \eqref{mi} follows, taking $A=(O-Q_2)/r_2$ and $B=(P_2-Q_2)/r_2$.

The point $P_1'$ lies on the line segment $OP_2$, and (recalling
that $\Ed(O,P_1)\ge ak^{1/3}$) is at Euclidean distance at least $c/2$ from both endpoints
of $OP_2$.
Thus, $P_1'$ lies on $XM$ where $X=O$ or $X=P_2$,
and $\Ed(X,P_1')\ge c/2$, so $\lambda=\Ed(X,P_1')/\Ed(X,M)\ge c_2k^{-1/3}$ for some
constant $c_2$.
From convexity of $d$ we have
\begin{eqnarray*}
 d(P_1',Q_2) &=& d\bigl((1-\lambda)X+\lambda M,Q_2\bigr) \\
 &\le& (1-\lambda) d(X,Q_2) + \lambda d(M,Q_2) \\
 &\le& (1-\lambda)r_2 + \lambda(1-\eps)r_2 = r_2-\eps\la r_2 = r_2-\Theta(1).
\end{eqnarray*}
Since $d(P_1,P_1')=o(1)$, it follows if $k$ is large enough that $d(P_1,Q_2)<r_2$,
contradicting our assumption that $d(P_i,Q_j)\ge r_j$ for all $i$ and $j$.
\end{proof}

The proof above can be easily adapted to any norm on $\RR^3$ by redefining the event $B_4$:
we must exclude a set of directions that form a neighbourhood of the set
of directions in which the unit sphere contains a line segment.
One can check that this latter set has measure zero for any norm, so it has neighbourhoods with
arbitrarily small measure, and we can ensure that $\Pr(B_4)=o(e^{-Ck})$ as above.

\medskip
We now turn to the proof of Theorem~\ref{th_couplerobust}. The key observation
is that a point of the $s$ by $s$ torus $\TT(s)$ is robustly black with respect to $(\PP_2^+,\PP_2^-)$
if it remains black when the points of $\PP_2^+$ are shifted away from the plane by a suitable
distance.
As a shifted Poisson process is again a Poisson process (but on a different set),
this allows us to construct the coupling so that `defects' arise essentially independently --
this is the reason why Theorem~\ref{th_couplerobust} is easier to prove than
the corresponding result for random 
Voronoi tessellations in the plane, Theorem 6.1 of~\cite{Voronoi}.
More precisely, we can deduce Theorem~\ref{th_couplerobust} from the ungainly
Lemma~\ref{l_c2} below. In this result, the notation is rather unnatural --
we write $\PP_1$ and $\PR_2$ for coloured Poisson processes where the probability that a point
is black is $p_1$ and $p_2$, respectively. We write $\PR_2$ rather than $\PP_2$
since we shall modify the process $\PR_2$ to obtain a process $\PP_2$
with the properties required for Theorem~\ref{th_couplerobust}.

\begin{lemma}\label{l_c2}
Let $0<p_1<p_2<1$ and $\eps'>0$ be given. Set $\delta'=s^{-\eps'}$.
We may construct in
the same probability space Poisson processes $\PP_1^+$, $\PP_1^-$,
and $\PR_2^-$ on $\TT(s)\times [0,s]$ of intensities $p_1$, $1-p_1$, and $1-p_2$,
respectively, and a Poisson process $\PR_2^+$ on $\TT(s)\times [\delta',s]$ of
intensity $p_2$, so that $\PP_1^+$ and $\PP_1^-$ are independent,
$\PR_2^+$ and $\PR_2^-$ are independent,
and the following global event $\Eglp$ holds \whp\ as $s\to\infty$: for
every piecewise-linear path $P_1\subset \TT(s)$ which is black with respect to
$(\PP_1^+,\PP_1^-)$ there is a piecewise-linear path $P_2\subset \TT(s)$ which is
black with respect to $(\PR_2^+,\PR_2^-)$, such
that every point of $P_2$ is within distance $\log s$ of some
point of $P_1$ and {\em vice versa}.
\end{lemma}

Before proving Lemma~\ref{l_c2}, let us show that it implies
Theorem~\ref{th_couplerobust}.  From here on we follow the convention
in~\cite{Voronoi} of writing, $x$, $x_i$ etc for points of $\TT(s)$,
and $z$, $z'$ etc for points of our Poisson processes in
$\TT(s)\times[0,s]$.  (The more natural notation $P\in \PP$ used above
becomes confusing when paths $P$ are involved.)  With this convention
a point $x$ is black if the nearest $z$ is black.

\begin{proof}[Proof of Theorem~\ref{th_couplerobust}]
We may assume without loss of generality that $\delta=s^{-\eps}<1$.
Set $\eps'=\eps/3$, and note that $\delta'=s^{-\eps/3}$ is larger than $\delta$.
Let $\PP_1^+,\PP_1^-,\PR_2^+$ and $\PR_2^-$ be coupled Poisson process with the properties
described in Lemma~\ref{l_c2}, set $\PP_2^-=\PR_2^-$, and let
\[
 \PP_2^+ = \{(x,y,z-\delta'): (x,y,z)\in \PR_2^+\} \cup \PP',
\]
where $\PP'$ is a Poisson process of intensity $p_2$ on $\TT(s)\times [s-\delta',s]$
that is independent of $\PR_2^\pm$.
Note that $\PP_2^+$ and $\PP_2^-$ are independent Poisson processes on $\TT(s)\times [0,s]$
with the desired intensities $p_2$ and $1-p_2$.
It remains only to show that $\Egl$ holds \whp.

From Lemma~\ref{l_c2}, we may assume that $\Eglp$ holds.
Let $E$ be the event that every point of $\TT(s)$
is within distance $O((\log s)^{1/3})=\OS(1)$ of some point of $\PR_2^-$, say.
(As usual, we write $f(s)=\OS(g(s))$ if $f(s)=O((\log s)^C g(s))$ for some constant $C$.)
Then, from basic properties of Poisson processes, $E$ holds \whp.
To complete the proof, we shall show that $\Eglp\cap E$ implies $\Egl$. 

Suppose that $\Eglp$ and $E$ hold, 
and let $P_1$ be a path that is black with respect to $(\PP_1^+,\PP_1^-)$.
Then, since $\Eglp$ holds, there is a path $P_2$ within Hausdorff distance $\log s$
of $P_1$ that is black with respect to $(\PR_2^+,\PR_2^-)$.
Let $x$ be any point of $P_2$, and let $z^+$ and $z^-$ be points of $\PR_2^+$ and $\PR_2^-$
at minimal distance from $x$. Since $x$ is black and $E$ holds, we have $d(x,z^+)\le d(x,z^-) = \OS(1)$.
Let $z'$ be the point of $\PP_2^+$ obtained by shifting $z^+$ by a distance $\delta'$
in the negative $z$-direction. We claim that
\begin{equation}\label{closer}
 d(x,z') \le d(x,z^+) - (\delta')^2/\OS(1).
\end{equation}
If $d=\JMd$, this is immediate; the reduction in distance is exactly $\delta'$.
If $d=\Ed$, then the extreme case is when $z'$ lies in $\TT(s)$ and $x$ and $z'$
are at maximal distance, in which case the claim follows from Pythagoras' Theorem.

Now $(\delta')^2=s^{-2\eps/3}$, while $\delta=s^{-\eps}$. If $s$ is large enough,
then as $\PP_2^-=\PR_2^-$, it follows from \eqref{closer} that
\[
 d(x,\PP_2^+) \le d(x,z') \le d(x,z^+)-2C_d\delta \le d(x,z^-)-2C_d\delta = d(x,\PP_2^-)-2C_d\delta.
\]
Thus $x$ is $(2C_d\delta)$-robustly black with respect to $(\PP_2^+,\PP_2^-)$.
Since $x$ was an arbitrary point of $P_2$, it follows that $P_2$ is $(2C_d\delta)$-robustly black.
Finally, as $P_1$ was arbitrary, $\Egl$ holds, as required.
\end{proof}

We now turn to the proof of Lemma~\ref{l_c2}.

\begin{proof}[Proof of Lemma~\ref{l_c2}]
Throughout the proof we write $\delta$ and $\eps$ for the quantities
$\delta'$ and $\eps'$ appearing in the statement of the lemma.

To construct our coupled Poisson processes,
we start with three independent Poisson processes,
a process $\PP$ of intensity $1$ on $\TT(s)\times[\delta,s]$,
and two processes on $\TT(s)\times[0,\delta]$,
a process $\PP_\delta$ of intensity $1$
and a process $\pDD$ of (much higher) intensity $\delta^{-1/2}$.
We shall form (homogeneous) Poisson processes
$\PP_1^\pm$ and $\PR_2^\pm$ with the properties described in the statement
of the lemma by assigning
every point of $\PP$ to exactly
one of $\PP_1^\pm$ and to exactly one of $\PR_2^\pm$,
assigning certain points of $\PP_\delta$ to $\PP_1^-$ and/or
$\PR_2^-$, and assigning certain points of $\pDD$ to $\PP_1^+$.
[Thinking of $\PP_1^\pm$ as a single Poisson process whose points are coloured
black/white according to a colouring $\col_1$, and writing $\PR_2^\pm$ similarly
in terms of ($\PR_2$, $\col_2$), then every point $z$ of $\PP$ is present in both
$\PP_1$ and $\PR_2$, although the colours $\col_1(z)$ and $\col_2(z)$ may be different.
Certain points of $\PP_\delta$ are present as white points in $\PP_1$ and/or $\PR_2$,
and certain points of $\pDD$ are black in $\PP_1^+$.]

Table~\ref{t1} summarizes the domains and intensities of these
Poisson processes, as well as two others that we shall consider in the proof.
The first four lines show the processes we shall construct; the remaining
lines concern processes used in the construction.

\begin{table}[htb]
\begin{center}
\begin{tabular}{|c|c|c|}
 \hline
Process & Domain & Intensity\\
 \hline
$\PP_1^+$ & $\TT(s)\times[0,s]$ & $p_1$ \\
$\PP_1^-$ & $\TT(s)\times[0,s]$ & $1-p_1$ \\
$\PR_2^+$ & $\TT(s)\times[\delta,s]$ & $p_2$ \\
$\PR_2^-$ & $\TT(s)\times[0,s]$ & $1-p_2$ \\
 \hline
$\PP$ & $\TT(s)\times[\delta,s]$ & 1 \\
$\PP_\delta$ & $\TT(s)\times[0,\delta]$ & 1 \\
$\pDD$ & $\TT(s)\times[0,\delta]$ & $\delta^{-1/2}$ \\
\hline
$\PPop$ & $\TT(s)\times[\delta,s]$ & $p_1$ \\
$\DD$ & $\TT(s)\times[0,\delta]$ & $p_1$ \\
\hline
\end{tabular}
\end{center}
\caption{The various Poisson processes involved in the statement and proof of Lemma~\ref{l_c2}.
We must construct $\PP_1^+$ and $\PP_1^-$ to be independent, and $\PR_2^+$ and
$\PR_2^-$ to be independent. The processes $\PP$, $\PP_\delta$ and $\pDD$
are independent by definition; $\PPop$ and $\PP_1^-$ will be independent by construction.}
\label{t1}
\end{table}

Table~\ref{t2} shows the probabilities with which the points of our independent
processes $\PP$, $\PP_\delta$ and $\pDD$ are included into the derived processes.
Note that $\PP_1^+$ is simply $\PPop\cup\DD$.

\begin{table}[htb]
\begin{center}
\begin{tabular}{|c|cc|cccc|}
 \hline
 & $\DD$ & $\PPop$ & $\PP_1^+$ & $\PP_1^-$ & $\PR_2^+$ & $\PR_2^-$ \\
 \hline
$\PP$ & & $p_1$ & $p_1$ & $1-p_1$ & $p_2$ & $1-p_2$ \\
$\PP_\delta$ & & & & $1-p_1$ & & $1-p_2$ \\
$\pDD$ & $p_1\delta^{1/2}$ & & $p_1\delta^{1/2}$ & & & \\
 \hline
\end{tabular}
\end{center}
\caption{The matrix of inclusion probabilities when $\DD$, $\PPop$, $\PP_1^\pm$ and $\PR_2^\pm$
are constructed from $\PP$, $\PP_\delta$ and $\pDD$. Note that $\PP_1^+=\PPop\cup\DD$.}
\label{t2}
\end{table}

In constructing $\PP_1^\pm$ and $\PR_2^\pm$ we shall define
two intermediate processes, $\PPop$ and $\DD$, whose union will form $\PP_1^+$.
The notation reflects the fact that $\PPop$ will consist of almost
all of $\PP_1^+$, more precisely, all points except those with $z$-coordinate
at most $\delta$. We write $\DD$ for the remaining points of $\PP_1^+$
since these points will form `defects' in our first attempt at a coupling.

By the `natural' coupling of $\PP_1^\pm$ with $\PR_2^\pm$
we shall mean the coupling obtained as follows.
Given $\PP$, $\PP_\delta$ and $\pDD$ as above,
for each point $z$ of $\PP$, toss a three-sided coin:
with probability $p_1$ assign $z$ to both
$\PPop$ and $\PR_2^+$, with probability $p_2-p_1$ assign $z$ to
$\PR_2^+$ and $\PP_1^-$, and with probability $1-p_2$ assign $z$ to
$\PR_2^-$ and $\PP_1^-$.
Similarly, for each point $z$ of $\PP_\delta$, with probability $p_1$
assign $z$ to neither $\PP_1^-$ nor $\PR_2^-$, with probability $p_2-p_1$
assign $z$ to $\PP_1^-$, and with probability $1-p_2$ assign $z$
to both $\PP_1^-$ and $\PR_2^-$.
Given
$\PP$ and $\PP_\delta$, we make all these choices independently.  Thus, the sets
$\PR_2^+$ and $\PR_2^-$ we obtain are independent Poisson processes on
$\TT(s)\times[\delta,s]$ and $\TT(s)\times[0,s]$, respectively, with
respective intensities $p_2$ and $1-p_2$.  Also, $\PPop$ and
$\PP_1^-$ are independent Poisson processes on
$\TT(s)\times[\delta,s]$ and $\TT(s)\times[0,s]$ respectively, with
respective intensities $p_1$ and $1-p_1$. To obtain $\PP_1^+$,
form a set $\DD$ by selecting each point of $\pDD$ independently
with probability $p_1\delta^{1/2}$, and set $\PP_1^+=\PPop\cup \DD$.

The construction above has the property that
$\PPop\subset \PR_2^+$ and $\PP_1^-\supset \PR_2^-$. Thus, any path $P_1$
in $\TT(s)$ that is black with respect
to $(\PPop,\PP_1^-)$ is black with respect to $(\PR_2^+,\PR_2^-)$.
Unfortunately, when we add the points of $\DD$ to $\PPop$ to obtain $\PP_1^+$,
this may
introduce new black paths with respect to $(\PP_1^+,\PP_1^-)$, which need
not be black with respect to $(\PR_2^+,\PR_2^-)$. For
this reason, we think of the points of $\DD$ as {\em defects}.
To deal with these defects, we shall need to adjust the coupling.

By a {\em potential defect} we mean a point of $\pDD$. Each potential
defect has only a small probability, $p_1\delta^{1/2}$, of becoming a real defect.
Even though there are many more potential defects than defects,
the density of potential defects
is still low: the projection of $\pDD$ onto $\TT(s)$ is a two-dimensional
Poisson process with intensity $\delta^{-1/2}\delta=\delta^{1/2}=s^{-\eps/2}$.

In constructing our final coupling, we shall condition on $\PP$, $\PP_\delta$ and $\pDD$.
Given these processes, we shall construct the remaining processes in a way
that respects the inclusion probabilities shown in Table~\ref{t2}:
the only changes we shall make to the `natural' coupling just defined
are to the coupling of inclusion choices associated to $\PP_1^\pm$ with inclusion
choices associated to $\PR_2^\pm$.  The first
step is to describe certain very unlikely `bad' events defined
in terms of $\PP$, $\PP_\delta$ and $\pDD$.
When one of these events holds, we shall complete the coupling arbitrarily (for example,
as above), and $\Eglp$ will not necessarily hold. As the bad events
will have probability $o(1)$, this will not be a problem.

Let $A$ be a large constant, and let $B_1$ be the event that there is some
point of $\TT(s)$ for which no point of $\PP$ lies
within $d$-distance $A(\log s)^{1/3}$. An elementary calculation shows
that $\Pr(B_1)=o(1)$ if $A$ is chosen large enough. (For $d=\Ed$, we may take $A=1$,
or indeed any constant such that $\pi A^3/3>1$.)

Let us say that two potential defects, i.e., points of $\pDD$, are {\em close}
if they are at $d$-distance at most $4A(\log s)^{1/3}$; this relation may
be taken to define
a graph on $\pDD$. By a {\em cluster} of potential
defects we mean a component of the resulting graph. Let $B_2$ be the event
that there is a cluster of potential defects containing more than $10/\eps$ potential
defects. If $B_2$ holds then, considering a connected set
of {\em exactly} $\ceil{10/\eps}$ potential defects, there is a $d$-disk
$D\subset \TT(s)$ of radius $r=\ceil{10/\eps}4A(\log s)^{1/3} = 
O((\log s)^{1/3})=\OS(1)$ containing the projections
of at least $\ceil{10/\eps}$ potential defects. But we may cover $\TT(s)$
with $O(s^2)$ $d$-disks $D_i$
of radius $2r$ so that any $d$-disk of radius $r$ is
contained in some $D_i$. The projection of $\pDD$ onto $\TT(s)$ is a Poisson process
of intensity $\delta^{1/2}=s^{-\eps/2}$, so the expected number of images in a given
$D_i$ is $\lambda=\delta^{1/2}\area(D_i)$, which is at most $s^{-\eps/3}$
if $s$ is large enough.
Hence,
\begin{eqnarray*}
 \Pr(B_2) &\le& \Pr\bb{\exists D_i\hbox{ containing }\ge \ceil{10/\eps}\hbox{ images} } \\
 &\le& O(s^2) \frac{\la^{\ceil{10/\eps}}}{\ceil{10/\eps}!}e^{-\la} = O\bb{s^2\la^{\ceil{10/\eps}}}
 = o(1).
\end{eqnarray*}

Let us say that a point $z'$ of $\pDD$ and a point $z$ of $\PP$ are
{\em potentially adjacent}
if the Voronoi cells of $z$ and $z'$ in the tessellation of $\TT(s)$ associated
to the point set $\PP\cup\{z'\}$ and metric $d$ meet.
Note for later that whenever the Voronoi cells of $z\in \PP$ and $z'\in\pDD$
defined with respect to some point set $X$ containing $\PP\cup\{z'\}$ meet,
then $z$ and $z'$ are potentially adjacent:
deleting points of $X$ to obtain $\PP\cup\{z'\}$ can only enlarge
the Voronoi cells associated to $z$ and $z'$. We shall apply this observation
later with $X=\DD\cup\PP\subset \PP_1^+\cup\PP_1^-$.
Note also that whether or not $z$ and $z'$ are potentially adjacent depends
only on $\PP$, $\PP_\delta$ and $\pDD$, not on $\PP_1^\pm$,
which we have not yet constructed.

Let $a$ be a small constant to be chosen later,
and let $B_3$ be the event that there is some $z'\in\pDD$ that is potentially
adjacent to at least $a\log s$ points $z_1,z_2,\ldots,z_k$ of $\PP$.

\begin{claim*}
For any choice of the constant $a>0$, we have $\Pr(B_3)=o(1)$ as $s\to\infty$.
\end{claim*}

Since $\Pr(B_1)=o(1)$, the claim follows if we show
that $\Pr(B_1^\cc\cap B_3)=o(1)$. We shall deduce this from Theorem~\ref{th_faces}.

For $z'\in \pDD$, let $E_{z'}$ be the event that there are $k=\ceil{a\log s}$
points $z_1,\ldots,z_k$ of $\PP$
and corresponding points $x_1,\ldots,x_k\in \TT(s)\times\{0\}$
such that $d(z',x_i)=d(z_i,x_i)=r_i\le A(\log s)^{1/3}$, and $d(z_j,x_i)\ge r_i$ for all $i$, $j$.
If $B_1^\cc\cap B_3$ holds, then so does $E_{z'}$ for some $z'\in \pDD$:
there is some $z'\in \pDD$ potentially adjacent to $k$ points $z_1,\ldots,z_k\in \PP$.
The Voronoi cells of $z'$ and $z_i$ defined with respect to $\PP\cup\{z'\}$
meet
at a point $x_i$ that is equidistant from $z'$ and $z_i$, with
$d(x_i,z)\ge d(x_i,z_i)=d(x_i,z')$ for every point $z\in \PP$, and in particular
for all $z_j$.
Since $B_1$ holds, some point of $\PP$ is within distance
$A(\log s)^{1/3}$ of $x_i$, so $r_i=d(x_i,z_i)\le A(\log s)^{1/3}$, and $E_{z'}$
holds.

Note for later that, if $B_1$ does not hold, any two potentially adjacent points
are within distance $2A(\log s)^{1/3}$.

To show that $\Pr(B_1^\cc\cap B_3)=o(1)$,
we shall condition on $\pDD$. Note that $\E(|\pDD|)=\delta^{-1/2}\delta s^2\le s^2/2$.
Since $|\pDD|$ has a Poisson distribution, it follows that
$|\pDD|\le s^2$ \whp; in proving the claim
we may thus assume that $\pDD$ is fixed and that $|\pDD|\le s^2$.
Let $z'$ be a point of $\pDD$. As $\PP$ and $\pDD$ are independent,
extending $\PP$ to a Poisson process on $\TT(s)\times \RR$ and translating
through the vector $-z'$, so that $z'$ is moved to the origin,
we may realize $\PP-z'$ as a subset of a Poisson process on $\TT(s)\times\RR$
of intensity $1$. In this Poisson process, if $E_{z'}$ holds,
then the shifted points $P_i=z_i-z'$ and $Q_i=x_i-z'$
have the properties described in Theorem~\ref{th_faces}, with $k=\ceil{a\log s}$.
Of course, this result
concerns a process on $\RR^3$ rather than on $\TT(s)\times \RR$. But the event
considered only involves points within distance $O((\log s)^{1/3})=o(s)$ of the origin,
so this makes no difference.
Hence, by Theorem~\ref{th_faces}, $\Pr(E_{z'})=o(e^{-3k/a})=o(s^{-3})$.
Considering the $O(s^2)$ points of $\pDD$ separately, it follows that $\Pr(B_1^\cc\cap B_3)=o(1)$.
As $\Pr(B_1)=o(1)$, this proves the claim.

\medskip
Let us say that a point $z'\in \pDD$ and a point $z\in \PP_\delta$ are {\em very close}
if they are within distance $2A(\log s)^{1/3}$.
Let $B_4$ be the event that some $z'\in \pDD$ is very close to at least $a\log s$
points $z\in \PP_\delta$. As the projection of $\PP_\delta$ onto $\TT(s)$
has intensity $\delta= s^{-\eps}$, it is easy to check that \whp\ no point
of $\pDD$ is very close to more than $10/\eps$
points of $\PP_\delta$. (The argument is similar to but simpler than that
for $\Pr(B_2)=o(1)$ above, so we omit the details.) It follows that $\Pr(B_4)=o(1)$.

\medskip
{F}rom now on we condition on $\PP$, $\PP_\delta$, and $\pDD$; we regard
these sets as fixed for the rest of the proof, and assume, as we may,
that none of $B_1,\ldots,B_4$ holds. To complete the construction of the coupling,
it remains to
assign each point $z$ of $\PP$ to $\PPop$ with probability $p_1$ and otherwise
to $\PP_1^-$, to assign $z$ to $\PR_2^+$ or $\PR_2^-$
with probabilities $p_2$ and $1-p_2$,
to assign each point $z$ of $\PP_\delta$ to $\PP_1^-$ with probability $1-p_1$
and to $\PR_2^-$ with probability $1-p_2$, to select points of $\pDD$ with probability
$p_1\delta^{1/2}$ to form $\DD$, and then to set $\PP_1^+=\PPop\cup\DD$.
As long as all choices involved in constructing $(\PP_1^-,\PP_1^+)$
are independent, $\PP_1^-$ and $\PP_1^+$ will have the right marginal distribution;
the same holds for $(\PR_2^-, \PR_2^+)$. If we can also ensure that $\Eglp$
always holds, the coupling will have the properties claimed
in the statement of the lemma.

Let $C_1,C_2,\ldots\subset \pDD$ be the clusters of potential defects. By assumption,
$B_2$ does not hold, so no $C_i$ contains more than $10/\eps$ points of $\pDD$.
Let $\Gamma(C_i)$ denote the set of points $z$ of $\PP\cup \PP_\delta$ with the property
that $z$ is potentially adjacent (if $z\in\PP$) or very close (if $z\in\PP_\delta$)
to one or more points $z'\in C_i$.
As $B_1$ does not hold, any pair of potentially adjacent points is at distance
at most $2A(\log s)^{1/3}$. Any pair of very close points is also separated by at most
this distance. Hence, if $z\in\Gamma(C_i)\cap \Gamma(C_j)$, then the clusters
$C_i$ and $C_j$ contain close points, contradicting the definition of a cluster.
Thus, the sets $\Gamma(C_i)$ are disjoint.
Also, as $B_3\cup B_4$ does not hold, we have $|\Gamma(C_i)|\le 2a\log s|C_i|\le 20a\eps^{-1}\log s$
for every $i$.

We shall construct our coupling independently for each set $C_i\cup \Gamma(C_i)$
in a manner described below. For each point $z\in \PP$ not in $\bigcup_i \Gamma(C_i)$,
we toss a three-sided coin, assigning $z$ to $\PPop$ and $\PR_2^+$ with probability
$p_1$, to $\PP_1^-$ and $\PR_2^+$ with probability $p_2-p_1$, and to $\PP_1^-$
and $\PR_2^-$ with probability $1-p_2$. These choices are independent for
different $z$, and also independent of the choices made within the $C_i\cup\Gamma(C_i)$;
this corresponds to the `natural' coupling described at the start of the proof.
Similarly, for each $z\in \PP_\delta$ not in $\bigcup_i\Gamma(C_i)$,
we include $z$ into both $\PP_1^-$ and $\PR_2^-$ with probability $1-p_2$,
and into $\PP_1^-$ only with probability $p_2-p_1$.

Let $C$ be a cluster of potential defects. Let $\ol{\DD}$, $\ol{\PPop}$, $\ol{\PP_1^\pm}$
and $\ol{\PR_2^\pm}$ be defined as in the natural coupling, but restricting
our attention to $C\cup\Gamma(C)$, i.e., starting from $C$, $\Gamma(C)\cap \PP_\delta$
and $\Gamma(C)\cap \PP$ in place
of $\pDD$, $\PP_\delta$ and $\PP$. Note that $\ol{\PP_1^\pm}$ and $\ol{\PR_2^\pm}$ have the correct distributions
for the restrictions of $\PP_1^\pm$ and $\PR_2^\pm$ to $C\cup\Gamma(C)$; however, to ensure that $\Eglp$
holds we shall have to adjust the coupling while keeping the marginal distributions fixed.

Let $B(C)$ be the `bad' event that $\ol{\DD}$ is non-empty, i.e., that
one of the potential defects in $C$ is an actual defect. Since $|C|\le 10/\eps$,
and each $z'\in C$ is included in $\ol{\DD}$ with probability $p_1\delta^{1/2}=p_1s^{-\eps/2}$,
we have $\Pr(B(C))\le s^{-\eps/3}$ for $s$ large.

Let $G(C)$ be the `good' event that every point of $\Gamma(C)$ is in $\ol{\PP_1^-}$
but no point of $\Gamma(C)$ is in $\ol{\PR_2^-}$. (Thus, each $z\in \Gamma(C)\cap\PP$
is in $\ol{\PR_2^+}$, while each $z\in \Gamma(C)\cap \PP_\delta$ does not appear
in $\ol{\PR_2^\pm}$.)
From the definition of the natural coupling, $G(C)$ has probability $(p_2-p_1)^{|\Gamma(C)|}$.
Recall that $|\Gamma(C)|\le 20a\eps^{-1}\log s$. By choosing the constant $a$ sufficiently small,
we may ensure that 
$(p_2-p_1)^{|\Gamma(C)|}\ge 2s^{-\eps/3}$. Hence, $\Pr(G(C))\ge 2\Pr(B(C))$,
and there is some event $G'(C)\subset G(C)\setminus B(C)$ with probability $\Pr(B(C))$.
[To be pedantic, we must modify our probability
space at this point. Having conditioned on $\PP$, $\PP_\delta$ and $\pDD$, and
restricting our attention to coordinates involving points in $C\cup\Gamma(C)$, we
are working in a finite probability space. To ensure that we can choose an event
$G'(C)$ with exactly the right probability, we should work in a space without atoms,
so we simply adjoin one extra random
variable $U_C$ for each cluster $C$, with $U_C$ uniform on $[0,1]$, say, and the new
variables independent of everything else.]

We define our final coupling by `crossing over' the natural coupling on the events
$B(C)$ and $G'(C)$: writing $\Omega$ for the probability space
on which we have defined $\ol{\PP_1^\pm}$ and $\ol{\PR_2^\pm}$,
let $f$ be a measure preserving bijection from $B(C)\cup G'(C)$ to itself,
mapping $B(C)$ into $G'(C)$ and {\em vice versa}.
We define the restrictions of $\PP_1^\pm$ to $C\cup\Gamma(C)$
to be equal to $\ol{\PP_1^\pm}$. 
Temporarily abusing notation by writing $\PR_2^\pm$ for the restriction of $\PR_2^\pm$
to $C\cup\Gamma(C)$,
for $\omega\notin B(C)\cup G'(C)$
we set $\PR_2^\pm(\omega)=\ol{\PR_2^\pm}(\omega)$, while for $\omega\in B(C)\cup G'(C)$
we set $\PR_2^\pm(\omega)=\ol{\PR_2^\pm}(f(\omega))$.
We extend the coupling to all clusters independently.
As the crossing over does not affect the marginal distribution of (the restriction to $C\cup \Gamma(C)$ of)
$\PR_2^\pm$, we end up with the correct marginal distributions
for $\PP_1^\pm$ and $\PR_2^\pm$. It remains to check that $\Eglp$ holds.

We claim that in our final coupling
\begin{equation}\label{mono}
 \PPop\subset \PR_2^+ \quad\hbox{ and }\quad\PP_1^-\supset \PR_2^-
\end{equation}
always holds. 
The `natural' coupling has this property, so we must check that
it is preserved by the local `crossing over' within $C\cup \Gamma(C)$.
Recall that $C\subset \pDD$ consists only of potential defects,
which never appear in any of the processes in \eqref{mono}: some
potential defects will be selected to form $\DD$, which will later be added
to $\PPop$ to form $\PP_1^+$. Recall also that $\Gamma(C)\subset \PP\cup\PP_\delta$.
Finally, recall that we only `cross over' at elements
$\omega\in B(C)\cup G'(C)$ of the state space.

Suppose first that $\omega\in G'(C)\subset G(C)\setminus B(C)$.
Then, since $\omega\in G(C)$, every point of $\Gamma(C)$ is in $\PP_1^-$,
and (hence) none is in $\PPop$. Thus the restriction of \eqref{mono} to $C\cup\Gamma(C)$ holds
in this case.

Suppose next that $\omega\in B(C)$. Then $f(\omega)\in G'(C)\subset G(C)$, so
every point of $\Gamma(C)\cap \PP$ is in $\PR_2^+$ and no point
of $\Gamma(C)$ is in $\PR_2^-$. Recalling that
points of $\PP_\delta$ never appear in $\PPop$, it again follows
that the restriction of \eqref{mono} to $C\cup \Gamma(C)$ holds. 
Since $C$ was an arbitrary cluster, this 
establishes \eqref{mono}.

\smallskip
From \eqref{mono}, if $x\in \TT(s)$ is black with respect to $(\PP_1^+,\PP_1^-)$
and the closest point of $\PP_1^+\cup\PP_1^-$ is not a defect (a point of $\DD=\PP_1^+\setminus \PPop$),
then $x$ is black with respect to $(\PR_2^+,\PR_2^-)$. It remains to show that we can 
modify a black path to avoid the defects;
we shall do this by taking a short detour around each defect.

For each defect $z'$, let $V_{z'}$ be the Voronoi cell of $z'$ in the tessellation $\VV$
of $\TT(s)$ defined with respect
to $\DD\cup \PP\subset \PP_1^-\cup\PP_1^+$. Note that $V_{z'}$ contains
the cell of $z'$ in the tessellation associated to $(\PP_1^-,\PP_1^+)$.
Let $\{z_1',\ldots,z_r'\}\subset \DD\subset \pDD$ be a maximal set of defects
such that the union $U$ of the cells
$V_{z_i'}$ is connected, and let $\ob{U}$ denote the external boundary of $U$.
As $B_1$ does not hold, no cell $V_{z'}$
has radius larger than $A(\log s)^{1/3}$, so the $z_i'$ lie in a single cluster $C$
of potential defects. Since $B_2$ does not hold, it follows  that $U$ has diameter at most $O((\log s)^{1/3})$.

\begin{claim*}
Every point of $\ob{U}$ is black with respect to $(\PR_2^+,\PR_2^-)$.
\end{claim*}

Before proving this claim, let us note that it implies that $\Eglp$ holds.
Let $P_1$ be a piecewise-linear path in $\TT(s)$ which
is black with respect to $(\PP_1^+,\PP_1^-)$.
Every point of $P_1$ not in one
of the sets $U$
is black with respect to $(\PR_2^+,\PR_2^-)$, since the nearest point of $\PP_1^+\cup\PP_1^-$ is not a defect.
On the other hand, every point of some set $U$ is surrounded by a black cycle,
$\ob{U}$, that lies within distance $\log s$.
Thus, whenever $P_1$ visits a set $U$, we may replace a part of $P_1$ by a part of $\ob{U}$, obtaining a
path $P_2$ that is black with respect to $(\PR_2^+,\PR_2^-)$, with every point of $P_1$ within distance $\log s$
of some point of $P_2$, and {\em vice versa}.
Shifting $P_2$ slightly, we may assume that $P_2$ is piecewise linear.
As $P_1$ was arbitrary, this shows that $\Eglp$ holds.

All that remains is to prove the claim.
Let $x\in \ob{U}$. Then $x$ is in some cell $V_z$ of $\VV$
with $z\notin\{z_1',\ldots,z_r'\}$. Since $U\cup V_z$ is connected,
the maximality of $\{z_1',\ldots,z_r'\}$ implies that the point $z$
is not a defect.
Therefore, recalling that the tessellation $\VV$ is defined with respect to $\DD\cup \PP$,
we have  $z\in \PP$. Furthermore, 
no point of $\PP$ is closer to $x$ than $z$ is.
In the tessellation $\VV$, the cell $V_z$ meets one of the cells $V_{z_i'}$, $z_i'\in C$.
Hence, $z$ is potentially adjacent to  $z_i'$, so $z\in \Gamma(C)$.
As the point $z_1'\in C$ is a defect, $B(C)$ holds, i.e., $\omega\in B(C)$.
Hence, $f(\omega)\in G(C)$. Since $z\in \Gamma(C)\cap \PP$,
it follows that $z\in \PR_2^+$. It remains to check that no point $z''$ of $\PR_2^-$
is closer to $x$ than $z$ is. As $\PR_2^-\subset \PP\cup\PP_\delta$,
and $z$ is the closest point of $\PP$ to $x$,
it suffices to consider points $z''\in \PP_\delta$. But $B_1$ holds,
so $d(x,z_i')=d(x,z)=d(x,\PP)\le A(\log s)^{1/3}$,
so any $z''\in\PP_\delta$ with $d(x,z'')<d(x,z)$
is very close to $z_i'$, and so lies in $\Gamma(C)$ and hence, from the definition of $G(C)$, cannot
be in $\PR_2^-$.
Hence, every $x\in\ob{U}$ is indeed black with respect to $(\PR_2^+,\PR_2^-)$, completing the proof
of the claim, and hence of Lemma~\ref{l_c2}.
\end{proof}

\section{Extensions}

For simplicity, Theorem~\ref{th1} was stated and proved for tilings of
$\RR^2$ generated by two specific norms on $\RR^3$, namely $\JMd$ and
$\Ed$. In the context of growing crystals, it is natural to consider
certain other norms, for example the $\ell_1$-norm on $\RR^3$,
corresponding to the sum of $\ell_1$ on $\RR^2$ and time. The statement
and proof of Theorem~\ref{th1} adapt immediately to this setting: the only
change required is the definition of the event $B_4$ in the proof of Theorem~\ref{th_faces}.
In fact, the same comment applies to any norm of the form $\norm{\cdot}\oplus\ell_1$,
where $\norm{\cdot}$ is a norm on $\RR^2$ which has the symmetries of the square;
these symmetries are assumed when we apply the Russo--Seymour--Welsh type
result from~\cite{Voronoi}

Theorem~\ref{th1} and its proof also extend to two-dimensional slices of $d$-dimensional
(usual) Voronoi tessellations, this time with simple changes: for example,
the bound $O((\log s)^{1/3})$ on the distance
to the nearest point of $\PP$ must be replaced by $O((\log s)^{1/d})$,
and $\TT(s)\times [0,s]$ is replaced by $\TT(s)\times [0,s]^{d-2}$. When we `shift' points,
it suffices to change one coordinate, so the set $\TT(s)\times [\delta',s]$ appearing
in the statement of Lemma~\ref{l_c2} is replaced by $\TT(s)\times [\delta',s]\times [0,s]^{d-3}$.

\bigskip
\noindent
{\bf Acknowledgement.}
Much of this research was carried out during a visit of the authors to
the Institute for Mathematical Sciences, National University of
Singapore; we are grateful to the Institute for its support.

\end{document}